\documentclass[12pt,a4paper,reqno]{amsart}

\usepackage{amssymb,mathtools,mathrsfs}

\edef\restoreparindent{\parindent=\the\parindent\relax}
\usepackage{parskip}
\restoreparindent

\usepackage[marginratio=1:1]{geometry}

\makeatletter
\@namedef{subjclassname@2020}{\textup{2020} Mathematics Subject Classification}
\makeatother

\newtheorem{thm}{Theorem}
\newtheorem{lem}{Lemma}

\newtheorem{cor}{Corollary}
\newtheorem{conj}{Conjecture}

\newtheorem{Thm}{Theorem}

\newtheorem{Lem}{Lemma}

\theoremstyle{definition}
\newtheorem{rem}{Remark}
\newtheorem{defn}{Definition}

\newcommand{\IC}{{\mathbb C}}
\newcommand{\ID}{{\mathbb D}}

\newcommand{\IS}{{\mathcal S}}
\newcommand{\CC}{{\mathcal C}}
\newcommand{\IK}{{\mathcal K}}
\newcommand{\real}{{\operatorname{Re}\,}}

\newcommand{\eit}{{e^{i\theta}}}

\newcommand{\IF}{\mathscr{H}}
\newcommand{\IG}{\mathscr{G}}

\newcommand{\be}{\begin{equation}}
	\newcommand{\ee}{\end{equation}}

\newcommand{\blem}{\begin{lem}}
	\newcommand{\elem}{\end{lem}}

\newcommand{\bdefn}{\begin{defn}}
	\newcommand{\edefn}{\end{defn}}

\newcommand{\bthm}{\begin{thm}}
	\newcommand{\ethm}{\end{thm}}

\newcommand{\bcor}{\begin{cor}}
	\newcommand{\ecor}{\end{cor}}

\newcommand{\bconj}{\begin{conj}}
	\newcommand{\econj}{\end{conj}}

\newcommand{\brem}{\begin{rem}}
	\newcommand{\erem}{\end{rem}}

\newcommand{\bpf}{\begin{proof}}
	\newcommand{\epf}{\end{proof}}

\begin{document}
	
	\bibliographystyle{abbrv}
	
	\title[Hardy spaces of harmonic quasiconformal mappings and Baernstein's theorem]
	{Hardy spaces of harmonic quasiconformal mappings and Baernstein's theorem}

	\author{Suman Das}
	\address{Suman Das\vskip0.05cm Department of Mathematics with Computer Science, Guangdong Technion - Israel
		Institute of Technology, Shantou, Guangdong 515063, P. R. China.}
	\email{suman.das@gtiit.edu.cn}
	
	\author{Jie Huang}
	\address{Jie Huang \vskip0.05cm Department of Mathematics with Computer Science, Guangdong Technion - Israel
		Institute of Technology, Shantou, Guangdong 515063, P. R. China.}
	\email{jie.huang@gtiit.edu.cn}
	
	\author{Antti Rasila}
	\address{Antti Rasila \vskip0.05cm Department of Mathematics with Computer Science, Guangdong Technion - Israel
		Institute of Technology, Shantou, Guangdong 515063, P. R. China, \vskip0.01cm and \vskip0.01cm Department of Mathematics, Technion - Israel
		Institute of Technology, Haifa 3200003, Israel.}
	\email{antti.rasila@gtiit.edu.cn; antti.rasila@iki.fi}

	

	\begin{abstract}
		Let $\IS_H^0(K)$, $K\ge 1$, be the class of normalized $K$-quasiconformal harmonic mappings in the unit disk. We obtain Baernstein type extremal results for the analytic and co-analytic parts of functions in the geometric subclasses of $\IS_H^0(K)$. We then apply these results to obtain integral means estimates for the respective classes. Furthermore, we find the range of $p>0$ such that these geometric classes of harmonic quasiconformal mappings are contained in the Hardy space $h^p$, thereby refining some earlier results of Nowak.
	\end{abstract}

	\maketitle
	\pagestyle{myheadings}
	\markboth{S. Das, J. Huang, and A. Rasila}{Hardy spaces of harmonic quasiconformal mappings and Baernstein's theorem}
	
	
	
	\section{Introduction and preliminaries}\label{sec1}
	
	\subsection{Notations and background}
	Let $\ID\coloneqq\{z \in {\mathbb C}:\, |z|<1\}$ be the open unit disk in the complex plane $\IC$. If $f$ is analytic in $\ID$,
	the \textit{integral means} of $f$ are defined as $$M_p(r, f) \coloneqq \left( \frac{1}{2\pi}\int_{0}^{2\pi}\vert f(r e^{i\theta}) \vert ^p\, d\theta \right)^{1/p} \quad \text{for} \quad 0 < p < \infty,$$ and $$M_\infty(r,f) \coloneqq \sup_{ \vert z \vert=r} \vert f(z) \vert.$$
	The function $f$ belongs to the \textit{Hardy space} $H^p$ $(0 < p \le \infty)$ if $$\|f\|_p \coloneqq \sup_{0\le r<1} M_p(r,f)<\infty.$$ 
	Detailed surveys on Hardy spaces and integral means can be found in the book of Duren \cite{Duren}. 
	
	A real-valued function $u(x,y)$, twice continuously differentiable in $\ID$, is called \textit{harmonic} if it satisfies the Laplace equation $$\Delta u =\frac{\partial^2 u}{\partial x^2}+\frac{\partial^2 u}{\partial y^2}= 0 $$ in $\ID$. 
	A complex-valued function $f=u+iv$ is harmonic in $\ID$ if $u$ and $v$ are real-valued harmonic functions in $\ID$. Every such function has a unique canonical representation $f=h+\overline{g}$, where $h$ and $g$ are analytic in $\mathbb{D}$ with $g(0)=0$. Analogous to the $H^p$ spaces, the \textit{harmonic Hardy spaces} $h^p$ are the class of functions $f$, harmonic in $\ID$, such that $\|f\|_p<\infty$. We refer to the book of Duren \cite{Duren:Harmonic} for the theory of planar harmonic mappings, and to the monograph of Pavlovi\'c \cite{Pavbook} for a concise survey on $h^p$ spaces.
	
	The harmonic function $f=h+\overline{g}$ is locally univalent and sense-preserving in $\ID$ if and only if the Jacobian $J_f(z)=\vert h'(z)\vert^2-\vert g'(z)\vert^2$ is positive for every $z \in \ID$. Denote by $\IS_H$ the class of all sense-preserving univalent harmonic functions $f=h+\overline{g}$ in $\ID$ with the normalizations  $h(0)=g(0)=h'(0)-1=0$. The class $$\IS_H^0 \coloneqq\{f=h+\overline{g} \in \IS_H : g'(0)=0\}$$ is compact and in a one-to-one relation with $\IS_H$. Let $\IK_H$, $\CC_H$, $\IS_H^*$ and $\IK_H(\theta)$ ($0\le \theta < \pi $) be the subclasses of $\IS_H$ consisting of harmonic mappings onto convex, close-to-convex, starlike, and convex-in-one-direction regions, respectively. We denote by $\IK_H^0$, $\CC_H^0$, $\IS_H^{*0}$, and $\IK_H^0(\theta)$ the corresponding compact classes. These classes were investigated first by Clunie and Sheil-Small \cite{CS}, and later by many authors (see \cite{Duren:Harmonic}).
	
	
	Abu-Muhanna and Lyzzaik \cite{Abu_Lyzzaik} showed that every harmonic mapping $f \in \IS_H$ is of class $h^p$ for some $p>0$. Their results were later improved by Nowak \cite{Nowak}, and recently, by the first author and Sairam Kaliraj \cite{DK3}. Notably, Nowak found the sharp ranges that $\IK_H \subset h^p$ for $p<1/2$, and $\CC_H \subset h^p$ for $p<1/3$.
	
	\subsection{Harmonic $K$-quasiconformal mappings}
	For $K\ge 1$, a sense-preserving harmonic function $f$ is said to be $K$-\textit{quasiregular} if the \textit{dilatation} $$D_f\coloneqq \frac{|f_z|+|f_{\bar{z}}|}{|f_z|-|f_{\bar{z}}|} \le K$$ throughout $\ID$. If $f=h+\overline{g}$, this condition is equivalent to saying that its \textit{analytic dilatation} $\omega \coloneqq g'/h'$ satisfies the inequality $$|\omega(z)| \le k <1 \quad (z\in \ID),$$ where \be\label{eq1}k\coloneqq \frac{K-1}{K+1}.\ee We say that $f$ is $K$-\textit{quasiconformal} if $f$ is $K$-quasiregular and homeomorphic in $\ID$. For convenience, we simply write that $f$ is quasiconformal if it is $K$-quasiconformal for some $K\ge 1$. One can find the $H^p$-theory for quasiconformal mappings in the paper of Astala and Koskela \cite{ast_kos}. 
	
	Let us denote $$\IS_H^0(K)\coloneqq \left\{f\in \IS_H^0: f \text{ is } K \text{-quasiconformal}\right\}.$$
	In \cite{wang_rasila}, Wang et al. constructed a harmonic $K$-quasiconformal Koebe-type function and
    suggested, without proof, that the Taylor series coefficients of every function $f \in \IS_H^0(K)$ are dominated in modulus by the coefficients of this new function. This has been verified very recently by Li and Ponnusamy \cite{Li_Pon}, for quasiconformal mappings in the geometric classes $\CC_H^0$, $\IS_H^{*0}$, and $\IK_H^0(\theta)$. Moreover, they have constructed a harmonic quasiconformal half-plane mapping
    which has the largest coefficients (in modulus) among the quasiconformal mappings in $\IK_H^0$.
	
	In view of this emergence of novel extremal functions for harmonic quasiconformal mappings, it becomes of natural interest to study other well-known extremal problems for the geometric subclasses of $\IS_H^0(K)$.
	
	\subsection{Baernstein's theorem}
	Let $\mathcal{S} \subset \IS_H^0$ denote the class of univalent analytic functions $f$ in $\ID$ with $f(0)=f'(0)-1=0$. A celebrated result on the growth of univalent functions is Baernstein's discovery that in the class $\IS$, the Koebe function has the largest integral means.
	
	\begin{Thm}\label{Baernstein_Analytic}{\rm \cite{AB1}}
		If $f \in \mathcal{S}$ and $0<p<\infty$, then for $0<r<1$, $$M_p(r,f) \le M_p(r,\mathbf{k}),$$ where $\mathbf{k}(z)=z/(1-z)^2$ is the classical Koebe function.
	\end{Thm}
	
	Baernstein's theorem was later extended to derivatives by Leung \cite{Leung} and Brown \cite{Brown} for certain subclasses of $\mathcal{S}$, most notably for the class $\CC$ of close-to-convex functions.
	\begin{Thm}\label{Leung_Brown}\cite{Brown,Leung}
		For $f\in \CC$ and $0<p<\infty$, we have $$M_p(r,f')\le M_p(r,\mathbf{k}'), \quad 0<r<1.$$
	\end{Thm}
	We note that Theorem \ref{Leung_Brown} is not true for the whole class $\IS$ (see \cite[p. 229]{Duren-book1}). We further refer to \cite{AB2,Girela1,Girela2,Nowak2} for similar extremal results for analytic functions. Recently, the first author and Sairam Kaliraj \cite{DK2} obtained Baernstein type inequalities for univalent harmonic functions, which had been unexplored till that point.
	
	The purpose of this paper is to produce Baernstein type results for harmonic quasiconformal mappings in the geometric subclasses of $\IS_H^0$. The main results and their implications are presented in Section \ref{sec_results}. In Section \ref{sec0}, we recall some results from the literature that will be useful for our purpose. The proofs of the main theorems and corollaries are given in Section \ref{sec2}.
	
	\section{Main results and related discussions}\label{sec_results}
	In what follows, we always assume that $0<r<1$, to avoid repetition. Also, we use $K$ and $k$ interchangeably by virtue of the relation \eqref{eq1}.
	\bthm\label{thm1}
	Suppose $0<p<\infty$ and $f=h+\overline{g} \in \IK_H^0$ is a $K$-quasiconformal mapping. Then $$M_p(r,h') \le M_p(r,H_k),$$ and $$M_p(r,g') \le M_p(r, G_k),$$ 
	where $$H_k(z)\coloneqq \frac{1+z}{(1-z)^2(1-kz)} \quad
	\text{and}\quad  G_k(z)\coloneqq \frac{kz(1+z)}{(1-z)^2(1-kz)}.$$
	\ethm
	
	\bthm\label{thm2}
	Let $0<p<\infty$. If the harmonic $K$-quasiconformal mapping $f=h+\overline{g}$ belongs to any of the geometric subclasses $\CC_H^0$, $\IS_H^{*0}$, $\IK_H^0(\theta)$, then $$M_p(r,h') \le M_p(r,\IF_k),$$ and $$M_p(r,g') \le M_p(r,\IG_k),$$
		where $$\IF_k(z)\coloneqq \frac{(1+z)^2}{(1-z)^3(1-kz)} \quad
	\text{and}\quad  \IG_k(z)\coloneqq \frac{kz(1+z)^2}{(1-z)^3(1-kz)}.$$
	\ethm
	
	\brem\label{rem0}
	It is important to note that Baernstein's method, although elegant and intricate, cannot be applied to obtain results of the form $$M_p(r,h) \le M_p(r,H), \quad M_p(r,g) \le M_p(r,G),$$ or more generally, $$M_p(r,f) \le M_p(r,F),$$ where $f=h+\overline{g}$ is harmonic and $F=H+\overline{G}$ is some extremal function. Baernstein's proof is critically based on the fact that for any $f\in \IS$, the function $$u(z)\coloneqq \log\left|\frac{1}{f^{-1}(z)}\right|=-\log|f^{-1}(z)|$$ is subharmonic. Clearly, this notion does not work for harmonic functions. Therefore, in some sense, Theorems \ref{thm1} and \ref{thm2} (for $h',\,g'$ instead of $h,\,g$) are the closest to Baernstein type results that one could expect in the harmonic case.
	\erem
	
	Nevertheless, we can easily give upper bounds for $M_p(r,f)$ (in terms of a single extremal function for each $k$) by using Theorems \ref{thm1} and \ref{thm2}.
	
	\bcor\label{cor1}
	Suppose $p\ge 1$ and $f=h+\overline{g}$ is a harmonic $K$-quasiconformal mapping. If $f \in \IK_H^0$, then we have the upper bound $$M_p(r,f) \le (1+k) \int_{0}^r M_p(s,H_k)\, ds,$$ where $H_k(z)$ is as in Theorem \ref{thm1}.
	\ecor
	
	\bcor\label{cor2}
	Let $p\ge 1$ and $f=h+\overline{g}$ be a harmonic $K$-quasiconformal mapping. If $f$ belongs to any of the geometric subclasses $\CC_H^0$, $\IS_H^{*0}$, $\IK_H^0(\theta)$, then $$M_p(r,f) \le (1+k) \int_{0}^r M_p(s,\IF_k)\, ds,$$ where $\IF_k(z)$ is as in Theorem \ref{thm2}.
	\ecor
	
	The next theorem is a refinement of the results of Nowak mentioned earlier.
	
	\bthm\label{thm3}
	If $f \in \IK_H$ is a $K$-quasiconformal mapping, then $$f\in h^p \quad \text{for} \quad p<1.$$ Similarly, if $f$ is a $K$-quasiconformal mapping in $\CC_H$ (or $\IS_H^{*}$, $\IK_H(\theta)$), then $$f\in h^p \quad \text{for} \quad p<1/2.$$
	\ethm
	
	\brem
	Theorem \ref{thm3} refines Nowak's sharp results that $\IK_H \subset h^p$ for $p<1/2$, and $\CC_H \subset h^p$ for $p<1/3$. In other words, quasiconformal mappings in the geometric subclasses of $\IS_H$ have significantly restricted growth, as one might expect. Perhaps more interesting is the comparison with a very well-known result of Astala and Koskela \cite{ast_kos}. Theorem 3.2 of \cite{ast_kos} asserts that a $K$-quasiconformal mapping of the unit disk belongs to the quasiconformal Hardy space $\mathcal{H}^p$ for $p<1/(2K)$, and this range of $p$ is the best possible. However, $$1/(2K) \ll 1/2 \quad \text{as} \quad K \gg 1,$$ i.e., the ranges of $p$ obtained in Theorem \ref{thm3} are better. This suggests that the study of Hardy spaces of harmonic quasiconformal mappings may lead to new non-trivial results.
	\erem

    \subsection*{Open problems}
    \begin{enumerate}
        \item It seems that Theorems \ref{thm1} and \ref{thm2} are not sharp, in the sense that there are no functions $f=h+\bar{g}$ in the respective classes such that $h'$ and $g'$ coincide with the obtained extremal functions. Therefore, the sharpness of these theorems remain open.
        \item In view of Theorem \ref{thm3}, it is tempting to ask if every $f\in \IS_H^0(K)$ is indeed in $h^p$ for $p<1/2$. This is a potentially interesting question that may have further implications.
    \end{enumerate}

	\section{Necessary results from the literature}\label{sec0}
	The proofs of Theorems \ref{thm1} and \ref{thm2} are fundamentally based on the use of the \textit{Baernstein star-function}, as defined below.
	\bdefn \cite{AB1}
	For a real-valued function $g(x)$ integrable over $[-\pi,\pi]$, the Baernstein star-function is defined as $$g^*(\theta) = \sup_{\vert E\vert=2\theta} \int_{E} g(x)\, dx \quad (0 \le \theta \le \pi),$$ where $\vert E\vert$ is the Lebesgue measure of the set $E \subseteq [-\pi, \pi]$.
	\edefn
	The relevance of the star-function in the study of integral means is contained in the following result of
	Baernstein.
	
	\begin{Lem}\label{lemAB1}{\rm \cite{AB1}}
		For $g,h \in L^1[-\pi, \pi]$, the following statements are equivalent.
		\begin{enumerate}
			\item[(a)] For every convex nondecreasing function $\Phi$ on $(-\infty, \infty)$, $$\int_{-\pi}^{\pi} \Phi(g(x))\, dx \le \int_{-\pi}^{\pi} \Phi(h(x))\, dx.$$
			\item[(b)] $g^*(\theta) \le h^*(\theta), \quad 0 \le \theta \le \pi$.
		\end{enumerate}	
	\end{Lem}
	
	Here we note some useful properties of the Baernstein star-function, which are due to Leung \cite{Leung}.
	
	\begin{Lem}\label{Leung1}
		For $g,h \in L^1[-\pi, \pi]$, $$[g(\theta)+h(\theta)]^* \le g^*(\theta) + h^*(\theta).$$ Equality holds if $g$, $h$ are both symmetric in $[-\pi, \pi]$ and nonincreasing in $[0, \pi]$.	
	\end{Lem}
	
	\begin{Lem}\label{Leung2}
		If $g$, $h$ are subharmonic functions in $\ID$ and $g$ is subordinate to $h$, then for each $r$ in $(0,1)$, $$g^*(r\eit) \le h^*(r\eit), \quad 0 \le \theta \le \pi.$$	
	\end{Lem}
	
	\begin{Lem}\label{Leung3}
		If $p(z)=e^{i\beta}+p_1z+\dots$ is analytic and has positive real part in $\ID$, then $$\left(\log\vert p(r\eit)\vert\right)^* \le \left(\log\left\vert\frac{1+r\eit}{1-r\eit}\right\vert\right)^*, \quad 0 \le \theta \le \pi.$$	 
	\end{Lem}
	
	\brem\label{rem1} An important feature in the proof of Lemma \ref{Leung3} is that a rotation factor does not affect the star-function. This observation will be suitably deployed in our proofs.\erem
	
	Also, we shall make use of the following two results in the proof of Theorem \ref{thm3}. 
	
	\begin{Lem} \label{lem1}\cite{DK3}
		Let $0 < p < 1$. Suppose $f=h+\overline{g}$ is a locally univalent, sense-preserving harmonic function in $\ID$ with $f(0)=0$. Then
		$$
		\| f\|_p^p \le C \int_{0}^{1} (1-r)^{p-1} M_p^p(r,h')\, dr,
		$$
		where $C$ is a constant independent of $f$.
	\end{Lem}
	
	\begin{Lem}\label{lemD1}{\rm \cite[p. 65]{Duren}}
		For each $p>1$, $$\int_{0}^{2\pi} \frac{d\theta}{|1 - r\eit|^{p}}  = O \left(\frac{1}{(1-r)^{p-1}}\right) \quad \text{as } r \to 1.$$
	\end{Lem}

	\section{Proofs}\label{sec2}
	\subsection{Proof of Theorem \ref{thm1}}
	Suppose $f=h+\overline{g} \in \IK_H^0$ is a $K$-quasiconformal mapping. It follows from Lemma 5.11 of \cite{CS} that there exist real numbers $\alpha$, $\beta$ such that $$\real \left\{ \left(e^{i\alpha}h'(z)+e^{-i\alpha}g'(z)\right)\left(e^{i\beta}-e^{-i\beta}z^2\right)\right\} > 0$$ for all $z \in \ID$. Let us write $$P(z) = \left(e^{i\alpha}h'(z)+e^{-i\alpha}g'(z)\right)\left(e^{i\beta}-e^{-i\beta}z^2\right).$$ Then $\real P(z) >0$ and $\vert P(0)\vert=1$. We see that $$\log\vert h'(z)\vert=\log\vert P(z)\vert+\log\left\vert\frac{1}{1+e^{-2i\alpha}w(z)}\right\vert+\log\left\vert\frac{1}{1-e^{-2i\beta}z^2}\right\vert,$$ and $$\log\vert g'(z)\vert=\log\vert h'(z)\vert+\log\vert \omega(z)\vert.$$
	
	The dilatation $\omega(z) = g'(z)/h'(z)$ satisfies $\omega(0)=0$ and $\vert \omega(z)\vert\le k$ for all $z \in \ID$. Therefore, appealing to Schwarz lemma (for the function $\omega(z)/k$), we find that $|\omega(z)|\le k|z|$ for every $z\in \ID$, and $|\omega'(0)|\le k$. It follows that $$\frac{1}{1+e^{-2i\alpha}\omega(z)}\prec \frac{1}{1-e^{-2i\alpha}kz},$$ where $\prec$ denotes the usual subordination. In view of Lemmas \ref{Leung1}--\ref{Leung3} and Remark \ref{rem1}, we have for $z=r\eit$ $(0 \le \theta \le \pi)$,
	\begin{align*}
		\left(\log\vert h'(z)\vert\right)^* & \le \left(\log\vert P(z)\vert\right)^*+\left(\log\left\vert\frac{1}{1+e^{-2i\alpha}\omega(z)}\right\vert\right)^*+\left(\log\left\vert\frac{1}{1-e^{-2i\beta}z^2}\right\vert\right)^*\\ & \le \left(\log\left\vert\frac{1+z}{1-z}\right\vert\right)^*+\left(\log\left\vert\frac{1}{1-e^{-2i\alpha}kz}\right\vert\right)^*+\left(\log\left\vert\frac{1}{1-e^{-2i\beta}z^2}\right\vert\right)^*\\ & = \left(\log\left\vert\frac{1+z}{1-z}\right\vert\right)^*+\left(\log\left\vert\frac{1}{1-kz}\right\vert\right)^*+\left(\log\left\vert\frac{1}{1-z^2}\right\vert\right)^*.
			\end{align*}
Since $\frac{1}{1-z^2} \prec \frac{1}{1-z}$, it follows from Lemma \ref{Leung2} that $$\left(\log\left\vert\frac{1}{1-z^2}\right\vert\right)^*\leq \left(\log\left\vert\frac{1}{1-z}\right\vert\right)^*.$$ Now, we observe that for fixed $r$, the functions $$\left|\frac{1+re^{i\theta}}{1-re^{i\theta}}\right|, \, \left|\frac{1}{1-kre^{i\theta}}\right|, \, \left|\frac{1}{1-re^{i\theta}}\right|$$ are symmetric in $[-\pi, \pi]$ and nonincreasing in $[0, \pi]$, so that the equality in Lemma \ref{Leung1} holds, and therefore
	\begin{align*}
	\left(\log\vert h'(z)\vert\right)^* & \leq \left(\log\left\vert\frac{1+z}{1-z} \cdot \frac{1}{1-kz} \cdot \frac{1}{1-z} \right\vert\right)^*\\ & = \left(\log\left\vert \frac{1+z}{(1-z)^2(1-kz)}\right\vert\right)^*=\left(\log|H_k(z)|\right)^*.
\end{align*}
	Similarly, \begin{align*}
		\left(\log\vert g'(z)\vert\right)^* &\le \left(\log\vert h'(z)\vert\right)^*+\left(\log\vert \omega(z)\vert\right)^*\\ & \le \left(\log\vert H_k(z)\vert\right)^*+\left(\log\vert kz \vert\right)^*\\ & = \left(\log\left\vert \frac{kz(1+z)}{(1-z)^2(1-kz)}\right\vert\right)^*=\left(\log|G_k(z)|\right)^*.
	\end{align*}
	Therefore, choosing $\Phi(x)=e^{px}$ in Lemma \ref{lemAB1}, we find that $$M_p(r,h') \le M_p(r,H_k),$$ and $$M_p(r,g') \le M_p(r,G_k).$$ This completes the proof.\qed
	
	\subsection{Proof of Theorem \ref{thm2}}
	Since $\IS_H^{*0},\, \IK_H^0(\theta) \subset \CC_H^0$, it is enough to prove the result for the class $\CC_H^0$. Let $f \in \CC_H^0$ be a $K$-quasiconformal mapping. Lemma 4 of \cite{WLZ} implies that there exist real numbers $\mu$, $\theta_0$ and an analytic function $Q(z)$ with positive real part, such that $$\real\left\{Q(z)\left[ie^{i\theta_0} \left(1-z^2\right)\left(e^{-i\mu}h'(e^{i\theta_0}z)+e^{i\mu}g'(e^{i\theta_0}z)\right)\right]\right\} > 0, \quad z\in \ID.$$ Let us write $$R(z)=Q(z)\left[ie^{i\theta_0} \left(1-z^2\right)\left(e^{-i\mu}h'(e^{i\theta_0}z)+e^{i\mu}g'(e^{i\theta_0}z)\right)\right].$$ It follows that
	$$
	\log\vert h'(e^{i\theta_0}z)\vert =\log \vert R(z)\vert  +\log\left\vert\frac{1}{Q(z)}\right\vert  +\log\left\vert\frac{1}{1-z^2}\right\vert +\log\left\vert\frac{1}{1+e^{2i\mu}\omega(e^{i\theta_0}z)}\right\vert,
	$$
	and $$\log\vert g'(z)\vert=\log\vert h'(z)\vert+\log\vert \omega(z)\vert.$$
	
	As before, the dilatation $\omega(z) = g'(z)/h'(z)$ satisfies $\omega(0)=0$ and $\vert \omega(z)\vert\le k|z|$ for all $z \in \ID$. Without any loss of generality, we may assume $\vert Q(0)\vert=1$, so that $\vert R(0)\vert=1$. Since $Q(z)$ has positive real part, so does $1/Q(z)$. Therefore, by Lemmas \ref{Leung1}--\ref{Leung3} and Remark \ref{rem1}, we find that
	\begin{align*}
		(\log\vert h'(z)\vert)^*  \le \left(\log\left\vert\frac{1+z}{1-z}\right\vert\right)^* & + \left(\log\left\vert\frac{1+z}{1-z}\right\vert\right)^*\\ & +\left(\log\left\vert\frac{1}{1-z^2}\right\vert\right)^*+\left(\log\left\vert\frac{1}{1-kz}\right\vert\right)^*,\end{align*}
	where $z=r\eit$ $(0 \le \theta \le \pi)$. This implies, like in the previous proof, that
	$$(\log\vert h'(z)\vert)^* \le \left(\log\left\vert \frac{(1+z)^2}{(1-z)^3(1-kz)}\right\vert\right)^* = \left(\log\vert \IF_k(z) \vert\right)^*.$$
	On the other hand, \begin{align*}
		\left(\log\vert g'(z)\vert\right)^* &\le \left(\log\vert h'(z)\vert\right)^*+\left(\log\vert \omega(z)\vert\right)^*\\ & \le \left(\log\vert \IF_k(z)\vert\right)^*+\left(\log\vert kz \vert\right)^*\\ & = \left(\log\left\vert \frac{kz(1+z)^2}{(1-z)^3(1-kz)}\right\vert\right)^*=\left(\log|\IG_k(z)|\right)^*.
	\end{align*}
	The proof is therefore completed via an appeal to Lemma \ref{lemAB1} with $\Phi(x)=e^{px}$.\qed
	
	\subsection{Proof of Corollary \ref{cor1}}
	Since $f=h+\overline{g}$ and $p\ge 1$, Minkowski's inequality implies $$M_p(r,f) \le M_p(r,h)+M_p(r,g).$$ For fixed $\theta \in [0,2\pi)$, we can write $$h(r\eit) = \int_{0}^{r} h'(s\eit) \eit\, ds,$$ as $h(0)=0$. Therefore, an appeal to Minkowski's integral inequality gives $$M_p(r,h) \le \int_{0}^r M_p(s,h')\, ds.$$ Similarly, we have $$M_p(r,g) \le \int_{0}^r M_p(s,g')\, ds \le k \int_{0}^r M_p(s,h')\, ds,$$ since $|g'(z)| \le k|h'(z)|$. Therefore, it follows from Theorem \ref{thm1} that $$M_p(r,f) \le (1+k) \int_{0}^r M_p(s,h')\, ds \le (1+k) \int_{0}^r M_p(s, H_k)\, ds.$$ The proof of Corollary \ref{cor2} is similar, hence omitted.\qed
	
	\brem\label{rem2}
	The upper bound in Corollary \ref{cor1} (and Corollary \ref{cor2}) can be further simplified for $p>1$. Indeed, we see that \begin{align*}
		M_p^p(s,H_k) & =  \frac{1}{2\pi} \int_{0}^{2\pi} \frac{|1+s\eit|^{p}}{|1-s\eit|^{2p}|1-ks\eit|^p}\, d\theta\\ & \le \frac{2^p}{2\pi(1-k)^p} \int_{0}^{2\pi} \frac{d\theta}{|1-s\eit|^{2p}}\\ & \le \frac{1}{(1-k)^p}  \frac{C}{(1-s)^{2p-1}},
	\end{align*}
	by Lemma \ref{lemD1}, where $C$ is an absolute constant that need not be the same in the following occurrences. This implies $$M_p(s,H_k) \le \frac{1}{(1-k)}  \frac{C}{(1-s)^{2-1/p}},$$ and therefore, \begin{align*}
		M_p(r,f) & \le (1+k) \int_{0}^r M_p(s,H_k)\, ds \\ & \le \left(\frac{1+k}{1-k}\right) \int_{0}^r \frac{C}{(1-s)^{2-1/p}}\, ds\\ & \le K \left(\frac{p}{p-1}\right) \frac{C}{(1-r)^{1-1/p}},
	\end{align*}
	as $K=\dfrac{1+k}{1-k}$.
	\erem

	\subsection{Proof of Theorem \ref{thm3}}
	Suppose $f =h+\overline{g}\in \IK_H$ is a $K$-quasiconformal mapping. Then the function \be \label{new1}f_0=\frac{f-\overline{\alpha f}}{(1-|\alpha|^2)} \in \IK_H^0 \quad \left(\alpha=g'(0)\right)\ee is $K_0$-quasiconformal for some $K_0\ge 1$ (see \cite[p. 10]{Li_Pon}). We write $f_0=h_0+\overline{g_0}$ and \be\label{new2}k_0\coloneqq \frac{K_0-1}{K_0+1}.\ee Now, let $1/2<p<1$. Applying Theorem \ref{thm1} and Lemma \ref{lemD1}, as in Remark \ref{rem2} above, we find that $$M_p^p(r,h_0')\le M_p^p(r,H_{k_0}) \le  \frac{1}{(1-k_0)^p} \frac{C}{(1-r)^{2p-1}},$$ where $C$ is an absolute constant, not necessarily the same at each occurrence. It follows from Lemma \ref{lem1} that \begin{align*}
		\| f_0\|_p^p & \le C \int_{0}^{1} (1-r)^{p-1} M_p^p(r,H_{k_0})\, dr\\ & \le C \int_{0}^{1} \frac{dr}{(1-r)^{p}}.
	\end{align*}
	This integral is finite for $p<1$. Therefore, $f_0\in h^p$ for $p<1$. It follows that $f-\overline{\alpha f} \in h^p$, and one needs to observe the inequality $$\|f-\overline{\alpha f}\|_p^p \ge \|f\|_p^p-|\alpha|^p\|f\|_p^p=\left(1-|\alpha|^p\right)\|f\|_p^p$$ to conclude that $f\in h^p$ for $p<1$.
	
	For the second part of the theorem, assume that $f=h+\overline{g}$ is a $K$-quasiconformal mapping in $\CC_H$ (or $\IS_H^{*}$, $\IK_H(\theta)$). Then the function $f_0=h_0+\overline{g_0}$, defined as in \eqref{new1}, belongs to $\CC_H^0$ and is $K_0$-quasiconformal for some $K_0\ge 1$. Suppose $1/3<p<1$. Theorem \ref{thm1} and Lemma \ref{lemD1} show that $$M_p^p(r,h_0')\le M_p^p(r,\IF_{k_0}) \le  \frac{1}{(1-k_0)^p} \frac{C}{(1-r)^{3p-1}},$$ where $C$ is an absolute constant and $k_0$ is defined as in \eqref{new2}. Thus, it follows from Lemma \ref{lem1} that \begin{align*}
		\| f_0\|_p^p & \le C \int_{0}^{1} (1-r)^{p-1} M_p^p(r,\IF_{k_0})\, dr\\ & \le C \int_{0}^{1} \frac{dr}{(1-r)^{2p}}.
	\end{align*}
	The last integral is finite for $p<1/2$, and therefore, $f_0\in h^p$ for $p<1/2$. This readily implies that $f\in h^p$, $p<1/2$, and the proof is complete. \qed

    We conclude with the remark that the inclusion of $\IS_H^0(K)$ in the harmonic Bergman space $a^p$ have been studied very recently in \cite{DR2}.

	\subsection*{Acknowledgements} The authors thank the anonymous referee for useful suggestions which led to the refinements of Theorems 1 and 2. The research was partially supported by the Li~Ka~Shing Foundation STU-GTIIT Joint Research Grant (Grant no. 2024LKSFG06) and the NSF of Guangdong Province (Grant no. 2024A1515010467).
	
		\bibliography{references}

\begin{thebibliography}{10}

\bibitem{Abu_Lyzzaik}
Y.~Abu-Muhanna and A.~Lyzzaik.
\newblock The boundary behaviour of harmonic univalent maps.
\newblock {\em Pacific J. Math.}, 141(1):1--20, 1990.

\bibitem{ast_kos}
K.~Astala and P.~Koskela.
\newblock {$H^p$}-theory for quasiconformal mappings.
\newblock {\em Pure Appl. Math. Q.}, 7(1):19--50, 2011.

\bibitem{AB1}
A.~Baernstein, II.
\newblock Integral means, univalent functions and circular symmetrization.
\newblock {\em Acta Math.}, 133:139--169, 1974.

\bibitem{AB2}
A.~Baernstein, II.
\newblock Some sharp inequalities for conjugate functions.
\newblock {\em Indiana Univ. Math. J.}, 27(5):833--852, 1978.

\bibitem{Brown}
J.~E. Brown.
\newblock Derivatives of close-to-convex functions, integral means and bounded
  mean oscillation.
\newblock {\em Math. Z.}, 178(3):353--358, 1981.

\bibitem{CS}
J.~Clunie and T.~Sheil-Small.
\newblock Harmonic univalent functions.
\newblock {\em Ann. Acad. Sci. Fenn. Ser. A I Math.}, 9:3--25, 1984.

\bibitem{DR2}
S.~Das and A.~Rasila.
\newblock On harmonic quasiregular mappings in {B}ergman spaces.
\newblock {\em Potential Anal.}, accepted, 2025.

\bibitem{DK2}
S.~Das and A.~Sairam~Kaliraj.
\newblock Growth of harmonic mappings and {B}aernstein type inequalities.
\newblock {\em Potential Anal.}, 60(3):1121--1137, 2024.

\bibitem{DK3}
S.~Das and A.~Sairam~Kaliraj.
\newblock Integral mean estimates for univalent and locally univalent harmonic
  mappings.
\newblock {\em Canad. Math. Bull.}, 67(3):655--669, 2024.

\bibitem{Duren:Harmonic}
P.~Duren.
\newblock {\em Harmonic mappings in the plane}, volume 156 of {\em Cambridge
  Tracts in Mathematics}.
\newblock Cambridge University Press, Cambridge, 2004.

\bibitem{Duren}
P.~L. Duren.
\newblock {\em Theory of {$H^{p}$} spaces}.
\newblock Pure and Applied Mathematics, Vol. 38. Academic Press, New
  York-London, 1970.

\bibitem{Duren-book1}
P.~L. Duren.
\newblock {\em Univalent functions}, volume 259 of {\em Grundlehren der
  mathematischen Wissenschaften [Fundamental Principles of Mathematical
  Sciences]}.
\newblock Springer-Verlag, New York, 1983.

\bibitem{Girela1}
D.~Girela.
\newblock Integral means and {BMOA}-norms of logarithms of univalent functions.
\newblock {\em J. London Math. Soc. (2)}, 33(1):117--132, 1986.

\bibitem{Girela2}
D.~Girela.
\newblock Integral means, bounded mean oscillation, and {G}elfer functions.
\newblock {\em Proc. Amer. Math. Soc.}, 113(2):365--370, 1991.

\bibitem{Leung}
Y.~J. Leung.
\newblock Integral means of the derivatives of some univalent functions.
\newblock {\em Bull. London Math. Soc.}, 11(3):289--294, 1979.

\bibitem{Li_Pon}
P.~Li and S.~Ponnusamy.
\newblock On the coefficients estimate of {$K$}-quasiconformal harmonic
  mappings.
\newblock {\em arXiv:2504.08284}, 2025.

\bibitem{Nowak2}
M.~Nowak.
\newblock Some inequalities for {BMOA} functions.
\newblock {\em Complex Variables Theory Appl.}, 16(2-3):81--86, 1991.

\bibitem{Nowak}
M.~Nowak.
\newblock Integral means of univalent harmonic maps.
\newblock {\em Ann. Univ. Mariae Curie-Sk\l odowska Sect. A}, 50:155--162,
  1996.

\bibitem{Pavbook}
M.~Pavlovi\'{c}.
\newblock {\em Introduction to function spaces on the disk}, volume~20 of {\em
  Posebna Izdanja [Special Editions]}.
\newblock Matemati\v{c}ki Institut SANU, Belgrade, 2004.

\bibitem{WLZ}
X.-T. Wang, X.-Q. Liang, and Y.-L. Zhang.
\newblock Precise coefficient estimates for close-to-convex harmonic univalent
  mappings.
\newblock {\em J. Math. Anal. Appl.}, 263(2):501--509, 2001.

\bibitem{wang_rasila}
Z.-G. Wang, X.-Y. Wang, A.~Rasila, and J.-L. Qiu.
\newblock On a problem of {P}avlovi\'c involving harmonic quasiconformal
  mappings.
\newblock {\em arXiv:2405.19852}, 2024.

\end{thebibliography}
	
\end{document}